\documentclass[12pt]{amsart}
\usepackage{amssymb,multicol,latexsym,epsfig}
\usepackage[mathscr]{eucal}

\allowdisplaybreaks

\newcommand{\aA}{\mathfrak{A}}  
\newcommand{\aB}{{\mathfrak{A}}^{\rm ass}}  
\newcommand{\bA}{{\bar{\mathfrak{A}}}} 
\newcommand{\bB}{{\bar{\mathfrak{B}}}}

\newcommand{\scalar}[2]{\bigl<#1\,,\,#2\bigr>}

\newcommand{\theA}{\mathcal{A}}

\def\bra#1{\langle{#1}|}

\makeatletter
 \@addtoreset{equation}{section}
\def\@secnumfont{\bfseries}
\makeatother

\newcommand{\N}[1]{N\!=\!#1}

\newcommand{\half}{\frac{1}{2}}

\newcommand{\oC}{\mathbb{C}}
\newcommand{\oN}{\mathbb{N}}

\newcommand{\oR}{\mathbb{R}}
\newcommand{\oZ}{\mathbb{Z}}

\newtheorem{Thm}{Theorem}[section]
\newtheorem{Cor}[Thm]{Corollary}
\newtheorem{Lemma}[Thm]{Lemma}

\newtheorem{Prop}[Thm]{Proposition}

\theoremstyle{definition}

\newtheorem{Rem}[Thm]{Remark}

\newtheorem{Example}[Thm]{Example}

{\noindent{\hfill\mbox{\rule{.5em}{.5em}}\,}\par\medskip}

\def\cF{\mathcal{F}}
\def\cG{\mathcal{G}}
\def\cH{\mathcal{H}}

\def\cL{\mathcal{L}}

\def\cT{\mathcal{T}}

\newcommand{\genp}[1]{{\cG^+}}
\newcommand{\genm}[1]{{\cG^-}}
\newcommand{\genpm}[1]{{\cG^\pm}}
\newcommand{\genh}[2]{{\cH^{#2}}}
\newcommand{\rep}[2]{{\mathsf{I}_{#2}}}
\newcommand{\Brep}[2]{{\mathsf{J}_{#2}}}
\newcommand{\brep}[2]{{\bar{\mathsf{I}}_{#2}}}
\newcommand{\vect}[1]{{\mathbf{#1}}}
\def \NN {{\vect{N}}}
\def \nn {{\vect{n}}}
\def \ve {{\vect{e}}}
\def \vu {{\vect{u}}}
\newcommand{\alg}[2]{{{W}_{#2}}}
\newcommand{\diag}{{\mathsf{diag}}}

\newcommand{\bgenp}[1]{{\bar{\cG}^+}}
\newcommand{\bgenm}[1]{{\bar{\cG}^-}}
\newcommand{\bgenpm}[1]{{\bar{\cG}^\pm}}
\newcommand{\bgenh}[2]{{\bar{\cH}^{#2}}}

\newcommand{\one}{1\kern-4pt1}

\newcommand{\charsl}[2]{{\chi_{
{{\phantom{h}\kern-3pt #2}}}^{\phantom{y}\kern-3pt #1}
}}  
\newcommand{\charn}[2]{{\omega_{
{{\phantom{h}\kern-3pt #2}}}^{\phantom{y}\kern-3pt #1}
}} 
\newcommand{\scharn}[2]{{\varpi_{
{{\phantom{h}\kern-3pt #2}}}^{\phantom{y}\kern-3pt #1}
}} 
\newcommand{\bcharsl}[2]{{\overline{\chi}_{
{{\phantom{h}\kern-3pt #2}}}^{\phantom{y}\kern-3pt #1}
}}
\newcommand{\bcharn}[2]{{\overline{\omega}_{
{{\phantom{h}\kern-3pt #2}}}^{\phantom{y}\kern-3pt #1}
}}

\newcommand{\chrct}{{\rm ch\,}}

\newcommand{\ket}[1]{|#1\rangle}

\def\bar{\overline}

\newcommand{\spfn}[1]{\mathop{\mathsf{U}_{#1}}}

\newcommand{\commut}[2]{\left[#1,\,#2\right]}

\renewcommand{\atop}[2]{\genfrac{}{}{0pt}{}{#1}{#2}}

\newcommand{\qbin}[3]{\left[\atop{#1}{#2}\right]_{#3}}

\newcommand{\qsup}[3]{\left(\atop{#1}{#2}\right)_{#3}}
\newcommand{\LL}{\mathbf{L}}

\def\tensor{\otimes}
\renewcommand{\d}{\partial}

\begin{document}
\raggedbottom
\hfuzz=1pt
\vfuzz=1.6pt
\addtolength{\baselineskip}{4pt}


\title{
  Coinvariants for Lattice VOAs and 
  $q$-Supernomial Coefficients}

\author{B.~L.~Feigin}
\address{Landau Institute for Theoretical Physics, Russian
  Academy of Sciences}
\author{S.~A.~Loktev}
\address{Independent University of Moscow}
\author{I.~Yu.~Tipunin}
\address{Tamm Theory Division, Lebedev Physics Institute,
  Russian Academy of Sciences}

\begin{abstract}
We propose an alternative definition of $q$-supernomial
coefficients as characters of coinvariants for one dimensional 
lattice vertex operator algebras. This gives a new formula for
$q$-supernomial coefficients.
Along the way we prove that the spaces of coinvariants form
a bundle over the configuration space of complex points (the
configuration space includes diagonals).
\end{abstract}


\maketitle

\thispagestyle{empty}

\setcounter{tocdepth}{3}

\vspace*{-36pt}

\begin{center}
  \parbox{.95\textwidth}{
    \begin{multicols}{2}
      {\footnotesize
  \setcounter{tocdepth}{1} \tableofcontents}
    \end{multicols}
    }
\end{center}


\section{Introduction}
In this paper we study the spaces of conformal blocks for one
dimensional lattice vertex operator algebras.

Spaces of conformal blocks or modular functor (see \cite{S}) play an
important role in conformal field theory. They can be defined for a
wide class of vertex operator algebras and they parameterizes the
states of the corresponding field theory. To avoid cumbersome notion
let us describe our approach for affine algebras (which corresponds to 
the WZNW theory).

\def \gg {{\mathfrak g}}
\def \ag {\widehat{\gg}}

First, recall the definition of the space of conformal blocks.
Let $\gg$ be a simple Lie algebra, by $\ag$ denote the central extension
of the Lie algebra of $\gg$--valued Laurent polynomials. Fix a level
$k$, that is, identify the central element of $\ag$ with the positive
integer $k$.
As a geometric part of the data consider a
complex curve $C$ and a set of pairwise distinct 
points $p_1, \dots, p_n \in C$ with 
a local coordinate in a small vicinity of each point. Then define the
Lie algebra $\gg^{C\setminus\{p_1, \dots, p_{n}\}}$ of $\gg$-valued
meromorphic functions on $C$ with poles only at $p_i$. Clearly we have
the inclusion
$$ \gg^{C\setminus\{p_1, \dots, p_{n}\}} \hookrightarrow \ag \oplus
\dots \oplus \ag,$$
mapping any $\gg$-valued function to its Laurent expansions at the
points $p_1, \dots, p_n$. Another part of the data is a set of integrable
representations $L_1, \dots, L_n$ of $\ag$ at level $k$. Their
external tensor product is a representation of $\ag \oplus
\dots \oplus \ag$ and, therefore, of 
$\gg^{C\setminus\{p_1, \dots, p_n\}}$.
Then the
space of conformal block is defined as the space of coinvariants
$$\left< L_1, \dots, L_n\right>^C_{p_1, \dots, p_n} = \left( L_1 \boxtimes
\dots \boxtimes L_n \right)/ \gg^{C\setminus\{p_1, \dots, p_n\}},$$
i.e. the quotient of $L_1 \boxtimes
\dots \boxtimes L_n$ by the action of $\gg^{C\setminus\{p_1, \dots,
p_n\}}$.
It is well known that these spaces are finite-dimensional and form a
flat family over the module space of curves with punctures.

\def \nn {{\mathfrak n}}

In \cite{FFu} another construction is proposed. Let, as above, $C$ be
a curve, $p_1, \dots, p_n$ and $p$ be pairwise distinct points on
$C$. Then consider the subalgebra
$\gg^{C\setminus p}(p_1, \dots, p_n) \subset \gg^{C\setminus p}$ of 
functions with zeroes at $p_1, \dots, p_n$ and the subalgebra 
$\gg_\nn^{C\setminus p}(p_1, \dots, p_n) \subset \gg^{C\setminus p}$
of functions whose values at $p_1, \dots, p_n$ belongs to the
nilpotent Borel subalgebra $\nn\subset \gg$. Then for any integrable
representation $L$ at level $k$ we have the natural
isomorphisms (see \cite{FKLMM})
\begin{align}\label{eq:decomp1}
 L /\gg_\nn^{C\setminus p}(p_1, \dots, p_n) & \cong 
\bigoplus_{L_1,\dots, L_{n}} 
\left<L_1, \dots, L_n, L\right>^C_{p_1, \dots, p_n,p},\\
 L /\gg^{C\setminus p}(p_1, \dots, p_{n}) & \cong 
\bigoplus_{L_1,\dots, L_{n}} \pi(L_1)\otimes \dots \otimes
\pi(L_{n}) \otimes
\left<L_1, \dots, L_{n}, L\right>^C_{p_1, \dots, p_{n},p},
\label{eq:decomp2}
\end{align}
where sum goes over all the sets of integrable representations at
level $k$ and $\pi(L_i)$ denotes the finite-dimensional representation
of $\gg$ generated by the highest vector of $L_i$.

So the space of coinvariants of a single representation $L$ contains
the spaces of conformal blocks where $L$ appears.
An advantage of this construction is that these spaces of coinvariants
is well defined if some of the points $p_1, \dots, p_{n}$
coincides. In this case we just consider zeroes with
multiplicities. For multiple points we have not
isomorphisms~\eqref{eq:decomp1},~\eqref{eq:decomp2},
but it is easy to show that these
spaces form a coherent sheaf over the variety $S^{n} C$ of sets of 
points on $C$. 

First question is whether the dimension of the space of coinvariants
preserves when some points coincides. Indeed, this question is local
and it is enough to check this when $C$ is the Riemann sphere.
The answer is positive for $\gg
= sl_2$ (see \cite{FKLMM} for $L /\gg^{C\setminus p}(p_1, \dots, p_{n})$
and  \cite{FKLMM2} for $L /\gg_\nn^{C\setminus p}(p_1, \dots, p_{n})$)
and negative for $\gg = E_8$. Also there is the positive answer for
the similar question about the $(2,m)$-minimal models of the Virasoro
algebras (see \cite{FFr}). But the general case is far from been
studied. Here we prove the equality of dimensions in the case of
one dimensional lattice VOAs.

\def \HH {{\mathcal H}}

Another question is about the natural filtration on these
spaces. Namely, the affine Lie algebra $\gg$ and 
its integrable representations are graded by the action of the energy
element. But the subalgebras defined above are not homogeneous (except
the case of Riemann sphere and $p_1 = \dots = p_{n}$), hence the
spaces of coinvariants can be considered only as filtered
spaces. Nevertheless we consider the Hilbert polynomials $\HH(q)$ of
the  adjoint
graded spaces and interpret them as a $q$-analog of the correspondent
Verlinde numbers (which are dimensions of spaces of conformal blocks).
This object is interesting even in the case of the Riemann sphere
(and well studied only in this case).
In \cite{FL} there is a conjecture that expresses $\HH(q)$ in terms of
generalized $q$-supernomial coefficients.
For $(2,m)$-minimal models of the Virasoro
algebras the polynomials $\HH(q)$ take part in the Andrews-Gordon
identities and therefore can be decomposed into sums of
$q$-multinomial coefficients as well. In this paper we prove a formula (see
Theorem~\ref{th:char}, (ii)) which expresses
$\HH(q)$ for one-dimensional lattice VOAs as a sum of  
$q$-supernomial coefficients. Note that this formula can be
considered as another definition of $q$-supernomials.

These two question are connected. Namely, suppose that the dimension
preserves when the points coincides. Then in the case of the Riemann
sphere it is enough to find $\HH(q)$ for the graded subalgebra with 
 $p_1 = \dots = p_{n}$. It makes possible to degenerate the algebra
together with the subalgebra into an abelian one or a rather simple VOA
and obtain a fermionic type (in terms of \cite{KKMM})
formula for $\HH(q)$. Such a formula exists for $\gg = sl_2$
(see \cite{FKLMM}, \cite{FKLMM2}, \cite{FKLMM3}), for 
$(2,m)$-minimal models of the Virasoro algebra (see \cite{FFr}) and for
one-dimensional lattice VOAs (see
Theorem~\ref{th:char}, (i)). 

In more general cases ($W$-algebras, lattice VOAs, etc.) it is not
evident which subalgebras should be chosen instead of
$\gg^{C\setminus p}(p_1, \dots, p_{n})$ and $\gg_\nn^{C\setminus
p}(p_1, \dots, p_{n})$. Let us explain what we are doing in the case of
affine Lie algebras.
Introduce some ``elementary''
graded subalgebras $\gg_i \subset \ag$. As the subalgebras are graded,
they consist of $\gg$-valued functions with some local conditions at the zero
point. Then define the ``fused'' subalgebra $\gg^{C\setminus p}_{i_1,
\dots, i_{n}} (p_1, \dots, p_{n})$ of functions with pole only in
$p$ and local conditions in $p_i$ prescribed by the subalgebras
$\gg_i$. Then we have
\begin{equation}\label{eq:decompg}
L /\gg_{i_1, \dots, i_{n}}^{C\setminus p}(p_1, \dots, p_n)  \cong 
\bigoplus_{L_1,\dots, L_{n}} \left(L_1 / \gg_{i_1}\right) \otimes
\dots \otimes \left(L_{n} / \gg_{i_{n}}\right)\otimes 
\left<L_1, \dots, L_{n}, L\right>^C_{p_1, \dots, p_{n},p}.
\end{equation}
In particular, the subalgebra $\gg_\nn^{C\setminus
p}(p_1, \dots, p_n)$ is fused from $n$ copies of the affine nilpotent
Borel subalgebra $\gg_\nn \subset \ag$ and the subalgebra $\gg^{C\setminus
p}(p_1, \dots, p_n)$ is fused from $n$ copies of the affine nilpotent
parabolic subalgebra $\gg_0 \subset \ag$, so \eqref{eq:decomp1} and
\eqref{eq:decomp2} follows from \eqref{eq:decompg}.

The paper is organized as follows.

In Sec.~\ref{sec:lvoa} we recall some facts about lattice vertex
operator algebras (``VOA'') and their representations.
In Sec.~\ref{sec:ass} we introduce an associative algebra that acts
in each representation of the VOA. In this way representations of the
VOA can be obtained as induced representations of the associative
algebra.

In Sec.~\ref{sec:verlinde} we introduce spaces of coinvariants and
give basic statements about their dimensions.
In Sec.~\ref{sec:element-sub} elementary subalgebras are introduced and
the fused subalgebras are studied.

Sec.~\ref{sec:algebras} and~\ref{sec:combinator} are rather technical.
In Sec.~\ref{sec:algebras} the character of some auxiliary space of
coinvariants is calculated. In Sec.~\ref{sec:combinator} we show that
characters of these auxiliary spaces can be expressed in terms of 
$q$-supernomial coefficients. This allows us to prove the main results
of the paper in Sec.~\ref{sec:main-result}.

We use the standard notation for the $q$-binomial coefficients
\begin{equation}\label{q-binom}
  \qbin n m q =
  \begin{cases}
    \frac{(q)_n}{(q)_{n-m}(q)_m},&n\geq m\geq0,\\
    0,&\text{otherwise},
  \end{cases}\qquad\text{where}\quad
  (q)_n=(1-q)\dots(1-q^n).
\end{equation}

\section{Lattice VOAs}
\subsection{\label{sec:lvoa}}
Recall that a lattice Vertex Operator Algebra (VOA) is determined
by the following data: $\oR^n$, $\scalar\cdot\cdot$, $\Gamma$, where
$\scalar\cdot\cdot$ is a scalar product in~$\oR^n$,
$\Gamma$ is a lattice in~$\oR^n$ and~$\scalar\cdot\cdot$ takes
integer values on vectors from~$\Gamma$. Such an algebra is generated
by vertex operators~$V(\vect{x},z)$, $\vect{x}\in\Gamma$.

In this paper we concentrate at one dimensional lattice VOAs
with positively defined scalar product.
Let~$\alpha$ be the basis vector of~$\Gamma$
and~$\scalar\alpha\alpha=p$, $p\in\oN$. Denote this VOA by~$\aA_p$ 
for~$p=1,2,3,4,\dots$.
The algebra~$\aA_p$ is generated by~$\genp{p}(z)=V(\alpha,z)$ 
and~$\genm{p}(z)=V(-\alpha,z)$.
The generators~$\genpm{p}$ are bosonic 
($\genpm{p}(z)\genpm{p}(w)-\genpm{p}(w)\genpm{p}(z)=0$) 
whenever~$p$ is even and fermionic
($\genpm{p}(z)\genpm{p}(w)+\genpm{p}(w)\genpm{p}(z)=0$)
otherwise and satisfy the following relations
\begin{align}
  \bigl(\genp{p}(z)\bigr)^2&=0\,,\qquad\bigl(\d\genp{p}(z)\bigr)^2=0\,,\quad
  \dots\,,\quad\bigl(\d^{p-3}\genp{p}(z)\bigr)^2=0\,,\\
  \bigl(\genm{p}(z)\bigr)^2&=0\,,\qquad\bigl(\d\genm{p}(z)\bigr)^2=0\,,\quad
  \dots\,,\quad\bigl(\d^{p-3}\genm{p}(z)\bigr)^2=0\,,
\end{align}
whenever $p$ is even and
\begin{align}
  \d\genp{p}(z)\genp{p}(z)&=0\,,\qquad\d^3\genp{p}(z)\genp{p}(z)=0\,,\quad
  \dots\,,\quad\d^{p-2}\genp{p}(z)\genp{p}(z)=0\,,\\
  \d\genm{p}(z)\genm{p}(z)&=0\,,\qquad\d^3\genm{p}(z)\genm{p}(z)=0\,,\quad
  \dots\,,\quad\d^{p-2}\genm{p}(z)\genm{p}(z)=0\,,
\end{align}
whenever $p$ is odd.

Introduce currents~$\genh{p}{m}$
for~$0\leq m\leq p-2$ as the currents appear in the OPE
\begin{equation}\label{eq:OPE}
  \genp{p}(z)\genm{p}(w)=\frac{\frac{1}{p}}{(z-w)^{p}}
 +\frac{\genh{p}{0}}{(z-w)^{p-1}}
 +\frac{\genh{p}{1}}{(z-w)^{p-2}}+\,\dots\,
 +\frac{\genh{p}{p-2}}{z-w}+\mbox{\rm regular terms}.
\end{equation}
Currents $\genh{p}{0}$ and $\genh{p}{1}$ correspond to $u(1)$-current
and energy--momentum tensor respectively.
Choose Virasoro~$\cT$ so that~$\genp{p}$ and~$\genm{p}$ 
be primary fields with conformal
weights~$p-1$ and~$1$ respectively
\begin{equation}
  \cT=\genh{p}{1}-\frac{p-1}{2}\d\,\,\genh{p}{0}
\,.
\label{eq:emt}
\end{equation}

The mode decomposition corresponding the choice~\eqref{eq:emt}
is
\begin{align}
\genp{p}(z)&=\sum_{n\in\oZ}\genp{p}_n\,z^{-n-p+1}\,,&
\genm{p}(z)&=\sum_{n\in\oZ}\genm{p}_n\,z^{-n-1}\,,\\
\label{eq:modes-h}
\genh{p}{m}(z)&=\sum_{n\in\oZ}\genh{p}{m}_n\,z^{-n-m-1}\,,&
\cT(z)&=\sum_{n\in\oZ}\cL_n\,z^{-n-2}\,.
\end{align}
This gives the following commutation relations
\begin{equation}
  \commut{\cL_n}{\genp{p}_i}=((p-2)n-i)\genp{p}_{i+n}\,,\qquad
  \commut{\cL_n}{\genm{p}_i}=-i\genm{p}_{i+n}\,,\qquad
  \commut{\cL_0}{\genh{p}{m}_i}=-i\genh{p}{m}_i\,.
\end{equation}

The algebra~$\aA_p$ admits the following family of automorphisms
(spectral flow) $\spfn{\theta}:$
\begin{equation}\label{gen-U}
   \begin{aligned}
    \genp{p}_n&\mapsto\genp{p}_{n+\theta},&{}
    \genm{p}_n&\mapsto\genm{p}_{n-\theta},\\
    \genh{p}{0}_n&\mapsto\genh{p}{0}_n+\theta\frac{1}{p}
    \delta_{n,0},&\qquad{}
    \genh{p}{1}_n&\mapsto\genh{p}{1}_n+\theta\genh{p}{0}_n
    +\half\theta(\theta-1)\frac{1}{p}\delta_{n,0},\\
    \genh{p}{m}_n&\mapsto\genh{p}{m}_n+
    \theta\,\genh{p}{m-1}_n+\half\theta(\theta-1)\,\genh{p}{m-2}_n
    +\frac{1}{6}\theta(\theta-1)(\theta-2)\,\genh{p}{m-3}_n
    +\dots\kern-250pt\\
    &&\kern-10pt+\frac{1}{m!}\theta(\theta-1)(\theta-2)\dots(\theta-m)
    \frac{1}{p}\delta_{n,0},\kern-210pt
  \end{aligned}  
\end{equation}
where $\theta\in\oZ$. The automorphisms with $\theta\in p\oZ$
play the role of the affine part of the Weyl group.

Note that in the case~$p=1$ the algebra~$\aA_p$ can be considered
as the infinite--dimensional Clifford algebra, in the case~$p=2$
as $\widehat{sl_2}$ at level~$k=1$ and in the case~$p=3$ as
the~$\N2$-superconformal algebra with the central charge~$c=\frac{1}{3}$.

Recall some facts about representations of lattice VOAs.
The category of representations of any lattice VOA is
semisimple, that is, any representation decomposes into the direct sum
of irreducible ones. The irreducible representations are enumerated by
the quotient $\Gamma^\vee/ \Gamma$, where $\Gamma^\vee$ is the dual
lattice. Namely, the representation corresponding to the element $\vect y
\in \Gamma^\vee$ is the space of descendants of~$V(\vect{y},0)$, and
if two elements of $\Gamma^\vee$ differs by an element of $\Gamma$
then the corresponding representations are isomorphic. It is also
known that these representations form a minimal model with the
Verlinde algebra isomorphic to the group ring of $\Gamma^\vee/ \Gamma$.

Let us describe it in detail for the one-dimensional case.
Denote the representation corresponding to~$V(-\frac{r}{\alpha},0)$
by~$\rep p r$ with~$0\leq r\leq p-1$. This representation
contains the highest weight vector~$\ket{p;r}$ satisfying the following
conditions
\begin{multline}
  \genp{p}_i\ket{p;r}=\genm{p}_j\ket{p;r}=0\,,\qquad i\geq-p+r+2\,,\quad
                                                     j\geq-r\,,\\
   \genh{p}{0}_0\ket{p;r}=-\frac{n}{p}\ket{p;r}\,,\qquad
   \cL_0\ket{p;r}=\frac{n}{2p}(r-p+2)\ket{p;r}\,.
\end{multline}
Another description of~$\rep{p}{r}$ is as follows
\begin{equation}\label{eq:voc}
\rep{p}{r}=\bigoplus_{m\in\oZ}\rep{p}{r}(m)
\end{equation}
where~$\rep{p}{r}(m)$ is the Fock module over~$\genh{p}{0}$
with the highest weight vector~$\ket{p;r+pm}$.
In this realization~$\genpm{p}:\rep{p}{r}(m)\to\rep{p}{r}(m\pm1)$
act as vertex operators.
We call the representation~$\rep{p}{0}$ a {\em vacuum} representation.

Representations~$\rep p r$ with~$0\leq r\leq p-1$ compose a minimal
model and the Verlinde algebra (the fusion algebra) is
\begin{equation}\label{eq:Veralg}
  \rep{p}{r}\cdot\rep{p}{s}=\rep{p}{r+s\mod p}
\end{equation}
The vacuum representation $\rep{p}{0}$ is the unit in the Verlinde
algebra.
Note also that the module~$\rep p r$
is isomorphic to the contragradient to~$\rep{p}{p-r}$.
This equips the Verlinde algebra with the scalar product.

Automorphisms $\spfn{\theta}$ naturally act on the set of 
representations~$\rep{p}{r}$ as
\begin{equation}
  \spfn{\theta}\rep{p}{r}=\rep{p}{r+\theta\mod p}\,.
\end{equation}
In particular, the action of $\spfn{1}$ on the set of 
representations~$\rep{p}{r}$ coincides with the fusion with~$\rep{p}{1}$.

Each space~$\rep{p}{r}$ is bi-graded by operators~$\genh{p}{0}_0$ 
and~$\cL_0$. 
Characters of representations~$\rep{p}{r}$ corresponding to
this bi-grading are by definition
\begin{equation}
  \chrct\rep{p}{r}(q,z)={\rm tr}\,z^{\cH_0}\,q^{\cL_0}\,,
\end{equation}
where the trace is calculated over~${\rep{p}{r}}$
and~\eqref{eq:voc} implies
\begin{equation}\label{eq:char-I}
  \chrct\rep{p}{r}(q,z)=\frac{z^{-\frac{r}{p}}
           q^{\frac{r}{2p}(r-p+2)}}{(q)_\infty}
  \sum_{n\in\oZ}\,z^n\,q^{\frac{p}{2}(n^2+n)-n(r+1)}
\end{equation}
For any graded subspace an analogous series can be considered.
We will denote them by the symbol~$\chrct$ in front of the symbol
of the subspace.

\subsection{\label{sec:ass}}
Representations of~$\aA_p$ can be considered from a more 
algebraic viewpoint. Operators~$\genpm{p}_n$ act on each
representation~$\rep p r$ therefore the 
associative algebra generated by~$\genpm{p}_n$ acts on each~$\rep p r$.
Denote this algebra by~$\aB_p$ and let us describe its defining relations.

Note first of all that~$\aA_p$ contains a Virasoro algebra
hence the algebra of vector fields on the circle acts by 
derivatives (commutators with~$\cL_j$) on~$\aB_p$.
By~$\cF_\lambda$ denote the representation of the algebra of vector
fields in weight $\lambda$ tensor fields 
$\oC[z,z^{-1}](dz)^\lambda$.
Then the space $\{\genp{p}_i,\ i\in\oZ\}=W_+$ is isomorphic
to~$\cF_{2-p}$ and $\{\genm{p}_i,\ i\in\oZ\}=W_-$ is isomorphic
to~$\cF_0$. This isomorphism is given by the identification
$\genp{p}_i=z_1^{i-p+2}(dz_1)^{2-p}$, $\genm{p}_i=z_2^{i}(dz_2)^{0}$.

Generators~$\genpm{p}_i$ commute whenever~$p$ is even,
anticommute otherwise and in both cases satisfy quadratic 
relations
\begin{equation}\label{eq:defrel-1}
\sum_{\alpha+\beta=n}p(\alpha,\beta)\genpm{p}_\alpha\genpm{p}_\beta=0\,,
\qquad n\in\oZ\,,
\end{equation}
for any polynomial $p(\alpha,\beta)$ of degree
not bigger than~$p-2$.  

The relations between~$\genp{p}_i$ and~$\genm{p}_j$
are also quadratic.
Introduce the space $W_0$ of all commutators of~$\genp{p}_i$ and~$\genm{p}_j$
for even~$p$ and anticommutators for odd~$p$.  Let us describe $W_0$
as a quotient of the tensor product $W_+\tensor W_-$.

We identify~$W_+\tensor W_-$ with the space of Laurent
polynomials~$\oC[z_1,z_1^{-1},z_2,z_2^{-1}](dz_1)^{2-p}(dz_2)^0$.
Then introduce the
subspaces
$$W_0^{(j)}=
\oC[z_1,z_1^{-1},z_2,z_2^{-1}](z_1-z_2)^{p-2-j}(dz_1)^{2-p}(dz_2)^0,
\qquad j\in\oZ.$$
Note that the quotient~$W_0^{(j)}/W_0^{(j-1)}$
is isomorphic to~$\cF_{-j}$ as a module over the algebra of
vector fields.

The OPE~\eqref{eq:OPE} implies that the space $W_0^{(-2)}$ acts on
any~$\rep{p}{r}$ by zero 
and $W_0^{(-1)}$ acts by a scalar.
The scalar can be found 
as follows. We know that $W_0^{(-1)}/W_0^{(-2)}\simeq\cF_1$, but~$\cF_1$
is the space of $1$-forms, so there exists the functional
\begin{equation}\label{eq:defrel-2}
  \beta:W_0^{(-1)}\longrightarrow\cF_1
         \stackrel{\mbox{\rm res}}{\longrightarrow}\oC
\end{equation}
Let $K\subset W_0^{(-1)}$ be the kernel of the
mapping~$\beta$. Then $K$ is a subspace in~$W_+\tensor W_-$ 
and we have $W_0=\left(W_+\tensor W_-\right)/K$.

\begin{Thm}
The defining relations in~$\aB_p$ are~\eqref{eq:defrel-1}
and relations between~$W_+$ and~$W_-$ following from~\eqref{eq:defrel-2}. 
\end{Thm}
We don't use this Theorem anywhere in the paper, so let us omit the
proof. But this consideration gives us the following description of
$\genh{p}{i}_s$ in terms of Laurent polynomials
\begin{equation}
  \genh{p}{i}_s=z_1^{-p+2}z_2^{s+i}(z_1-z_2)^{p-i-2}(dz_1)^{p-2}(dz_2)^0
\end{equation}
This notion is very convenient for us. In particular, 
automorphism~\eqref{gen-U} acts on the space of Laurent polynomials as 
the multiplication with~$(z_1/z_2)^\theta$.

One can show that the category of representations of~$\aB_p$ with
highest weight is equivalent to the category of representations of~$\aA_p$.
So irreducible representations of~$\aB_p$ with highest weight
are precisely~$\rep{p}{r}$.

Define for the algebra~$\aB_p$ some class of representations. Each of these 
representations contains a cyclic vector that is 
annihilated by some subspaces in~$W_+\bigoplus W_0\bigoplus W_-$,
which we now describe.
Let~$\vect N =(N_+,N_-;N_0,N_1,\dots,N_{p-2})$ be a $p+1$-dimensional
vector with integer components that satisfy
\begin{equation}\label{eq:subalg}
  0\leq N_0\leq N_1\leq\dots\leq N_{p-2}\leq N_++N_-\,.
\end{equation}
We define the subspace $W(\vect N)=W_+(\vect N)\oplus W_0(\vect N)\oplus 
W_-(\vect N)$ as follows
\begin{align*}
  W_+(\vect N) = \,&\,\oC[z_1]z_1^{N_+}(dz_1)^{2-p}\subset
  \oC[z_1,z_1^{-1}](dz_1)^{2-p}\,,\\
  W_-(\vect N)  = \,&\,\oC[z_2]z_2^{N_-}(dz_2)^0\ \ \, \subset
  \oC[z_2,z_2^{-1}](dz_2)^0
\end{align*}
and the space~$W_0(\vect N)$ is the sum (not direct!) of the
images of the subspaces 
$$W_0^{(j)}(\vect N) =
\oC[z_1,z_2]z_1^{N_j-N_-}z_2^{N_-}(z_1-z_2)^{p-j-2}(dz_1)^{2-p} (dz_2)^0 
$$ 
for $j = 0\dots p-2$ 
under the mapping $W_+\tensor W_-\to W_0$. 

Note that if $N_0 = N_1 = \dots = N_{p-2} = N_+ + N_-$ then $W_0(\vect
N)$ is the subspace of (anti)commutators of $W_+(\vect N)$ and
$W_-(\vect N)$.

Clearly we have $\spfn{\theta}(W(\vect N))=W(\vect N+\theta\vect{u})$, where
\begin{equation}\label{eq:vect-u}
  \vect{u}=(1,-1;0,\dots,0)
\end{equation}

\begin{Prop}
  The irreducible representation $\rep{p}{r}$ is induced
  from the subspace~$W(r\vect u)$.
\end{Prop}

\begin{Example}
  Let $p=2$. Then we can identify $\genp{p}_i=e\cdot t^i$, $\genm{p}_i=f\cdot
  t^i$, $\genh{p}{0}_i=h\cdot t^i$, where $e,h,f$ is the standard
  basis of~$sl_2$. So our VOA is isomorphic to~$\widehat{sl_2}$ 
  at level~$1$. Then we have 
  \begin{equation}
    W((N_+,N_-,N_0))=e\cdot t^{N_+}\oC[t]\oplus h\cdot t^{N_0}\oC[t]
     \oplus f\cdot t^{N_-}\oC[t] 
  \end{equation}
 So for $k=1$ we cover the cases considered in \cite{FKLMM,FKLMM2,FKLMM3}.
\end{Example}

Let $R(\vect N)$ be the representation induced from the 
trivial representation of the subalgebra generated by the 
subspace~$W(\vect N)$.
We know that the representation $R(\vect N)$ is a direct sum of $\rep{p}{r}$.
The multiplicity of~$\rep{p}{r}$ in~$R(\vect N)$ is the dimension
of the space~${\rm Hom}_{\aB_p}(R(\vect N),\rep{p}{r})$, which
is isomorphic to the space of $W(\vect N)$-invariants in~$\rep{p}{r}$.

\section{Coinvariants and Verlinde rules}
\subsection{\label{sec:verlinde}}
It is more convenient to study the dual space to the space of
$W(\vect N)$-invariants in~$\rep{p}{r}$ namely the space of 
$W(\vect N)$-coinvariants in~$\rep{p}{r}^*$, which is by 
definition~$\rep{p}{r}^*/W(\vect N)\rep{p}{r}^*$. Note that the dual
representation $\rep{p}{r}^*$ differs from $\rep{p}{r}$ only by exchanging
positive and negatives modes.

\def \va {\vect a}

The space $W(\vect{N})$ can be deformed, i.e.{} can be
included in a family of subspaces of~$W$ depending on a set of parameters
that are points on~$\oC$. Given $k$ 
vectors~$\vect N^i=(N_+^i,N_-^i;N_0^i,N_1^i,\dots,N_{p-2}^i)$
satisfying~\eqref{eq:subalg}
and given a vector $\va = (a^1, \dots, a^k) \in \oC^k$
we define 
$W(\vect N^1,\dots, \vect N^k;\va)$ as the direct sum of the spaces
$W_{\pm}(\vect N^1,\dots, \vect N^k;\va)$ and $W_0(\vect N^1,\dots, \vect
N^k;\va)$, where
\begin{align*}
  W_+(\vect N^1,\dots, \vect N^k;\va)
=&\,\,
\oC[z_1](dz_1)^{2-p} \prod_{i=1}^k (z_1 - a_i)^{N_+^i}\subset
  \oC[z_1,z_1^{-1}](dz_1)^{2-p}\,,\\
 W_-(\vect N^1,\dots, \vect N^k;\va)
=&\,\,
\oC[z_2](dz_2)^0 \prod_{i=1}^k (z_2 - a_i)^{N_-^i}\subset
  \oC[z_2,z_2^{-1}](dz_2)^0
\end{align*}
and the space
$W_0(\vect N^1,\dots, \vect N^k;\va)$ is the sum (not direct!) of the
images of the subspaces 
$$W_0^{(j)}(\vect N^1,\dots, \vect N^k;\va) =
\oC[z_1,z_2](z_1-z_2)^{p-j-2}(dz_1)^{2-p} (dz_2)^0 \prod_{i=1}^k (z_1 -
a_i)^{N_j^i-N_-^i}(z_2 - a_i)^{N_-^i}$$
for $j = 0\dots p-2$ under the
mapping
$W_+\tensor W_-\to W_0$. 

Then $W(\vect N)=W(\vect N^1,\vect N^2,\dots, \vect N^k;0)$ for $\vect
N=\vect N^1+\dots +\vect N^k$, so the
family of
subalgebras 
$W(\vect N^1,\dots, \vect N^k;\va)$ is a deformation 
of~$W(\vect N)$.

Such a deformation can be applied to the calculation of dimensions for the
coinvariants. To simplify the notation denote by $\rep{p}{r}[W]$ 
the space of coinvariants
$\rep{p}{r}^*/ W\rep{p}{r}^*$.

\begin{Thm}\label{thm:Verlinde}
Suppose that points $a^1,a^2,\dots,a^k$ are pairwise distinct. Then the
dimension of $\rep{p}{r}[W(\vect
N^1,\dots,\vect N^k;\va)]$ coincides with the coefficient
before $\rep{p}{r}$ in the following expression
\begin{equation}\label{eq:FF}
\prod_{i=1}^k \left( \dim \rep{p}{0}[W(\vect
N^i)] \cdot \rep{p}{0}  + \dots + \dim \rep{p}{p-1}[W(\vect
N^i)] \cdot \rep{p}{p-1} \right),
\end{equation}
where product means the Verlinde multiplication.
\end{Thm}
\begin{proof}
First recall the fusion procedure in the case of one input and $k$
outputs.

\def \MM {\vect M}

Suppose that $\NN^1, \dots, \NN^k$ and $\MM^1, \dots, \MM^k$ are two sets
of vectors such that $\NN^i_j \ge \MM^i_j$ for any $i=1\dots k$,
$j=+,-,0,\dots, p-2$. Then we have the natural projection
$\rep{p}{r}[W(\vect
N^1,\dots, \vect N^k;\va)] \to \rep{p}{r}[W(\vect M^1,\dots, \vect
M^k;\va)]$.
Therefore the dual spaces form an inductive system and we can consider the
limit
$$\rep{p}{r}[\va] = \lim_\to \left(\rep{p}{r}[W(\vect 
N^1,\dots, \vect N^k;\va)]\right)^*.$$

The action of the algebra $\aB_p$ on $\rep{p}{r}$ can be extended to the
action  of the algebra $\aB_p \oplus \dots \oplus \aB_p$ on the limit
$\rep{p}{r}[\va]$, where summands correspond to points $a^1, \dots, a^k$.

The statement that the set of representations $\rep{p}{r}$ 
form a minimal model means that the representation~$\rep{p}{r}[\va]$
decomposes into the direct sum of external tensor products
$$\rep{p}{r}[\va] = \bigoplus_{\atop{0\le s_1, \dots, s_k <p}{\rep{p}{s_1}
\cdot \rep{p}{s_2}\cdot \dots \cdot \rep{p}{s_k} = \rep{p}{r}}}
\rep{p}{s_1} \boxtimes \dots
\boxtimes \rep{p}{s_k}$$
Emphasize that the sum goes over the set of representations whose
product in the Verlinde algebra is~$\rep{p}{r}$.

Clearly, the space $\rep{p}{r}[W(\vect
N^1,\dots,\vect N^k;\va)]$ is dual to the space of $W(\vect
N^1,\dots,\vect N^k;\va)$-invariants in $\rep{p}{r}[\va]$. As the
subspace $W(\vect
N^1,\dots,\vect N^k;\va)$ is dense with respect to the topology
of direct limit in the subspace $W(\NN^1) \oplus
\dots \oplus W(\NN^k)\subset \aB_p \oplus \dots \oplus \aB_p$, we have
$$\dim \rep{p}{r}[W(\vect
N^1,\dots,\vect N^k;\va)] = \sum_{\atop{0\le s_1, \dots, s_k
<p}{\rep{p}{s_1}\cdot \rep{p}{s_2}
\cdot \dots \cdot \rep{p}{s_k} = \rep{p}{r}}}
\dim \rep{p}{s_1}[W(\vect
N^1)] \cdot \dots \cdot \dim \rep{p}{s_k}[W(\vect N^k)],$$
that is equivalent to the statement of the theorem.
\end{proof}

Here we have a family of subalgebras parametrized by complex numbers $a^1,
\dots, a^k$. The deformation argument (see for example \cite{FKLMM2})
implies that dimension of the space of coinvariants with respect to a
special subalgebra is not less than one with respect to a generic
subalgebra.

\begin{Prop}\label{Prop:defarg}
The dimension of $\rep{p}{r}[W(\vect
N^1,\dots,\vect N^k;\va)]$ is minimal whenever $a^i$ are
pairwise distinct and maximal whenever $a^i$ are equal each to other.
\end{Prop}

\subsection{\label{sec:element-sub}}
Introduce elementary spaces $\alg{p}{i,j}=W(\vect N^{i,j})$ with 
$$\vect N^{i,j}=
(i-j,j;1,2,\dots,i-1,i,\dots,i)$$
and $0\leq i\leq p$, $j\in\oZ$.
\begin{Prop}\label{?}
We have
 \begin{equation}\label{eq:dims}
   \dim \rep{p}{r}[\alg{p}{i,j}] =
\sharp\{n | 0 \le pn + j +r \le i\}
 \end{equation}
\end{Prop}
In other words, 
$$\dim \rep{p}{r}[\alg{p}{i,j}] =\begin{cases}
2 & i=p,\  (r+j) \ {\rm mod}\ p = 0\\
1 & \mbox{otherwise}\\
0 & (r+j)\ {\rm mod} \ p \  > \  i
\end{cases}$$
\begin{proof}
  Note that $\alg{p}{i,j}$ contains the subalgebra
  $\genh{p}{0}_{\geq1}$
therefore the space of $\alg{p}{i,j}$-invariants in~$\rep p r$
is spanned by several vectors~$\ket{p;r+pn}$, $n\in\oZ$.
Taking the invariants of~$\genp{p}_{\geq i-j-p+2}$ and~$\genm{p}_{\geq j}$
selects the vectors~$\ket{p;r+pn}$ with~$-j\leq r+pn\leq i-j$.
We need to show that these vectors are invariant.

Note that if $-j\leq r+pn\leq i-j$ then
$\alg{p}{i,j}\subset\alg{p}{0,-(r+pn)}$. And we know that the vector
$\ket{p;s}$ is annihilated by the subspace $\alg{p}{0,-s}$.
\end{proof}

According to \eqref{eq:FF}, \eqref{eq:dims} and \eqref{eq:Veralg} we
define the positive integers~$d_{p,r}[i^1,j^1;i^2,j^2;\dots;i^k,j^k]$
by the following formula
\begin{equation}\label{eq:def-d}
  \prod_{s=1}^k(\xi^{-j^s} + \xi^{-j^s+1} + \dots + \xi^{-j^s + i^s}) 
  =\sum_{r=0}^{p-1}d_{p,r}[i^1,j^1;i^2,j^2;\dots;i^k,j^k]\cdot \xi^r,
\end{equation}
where $\xi$ is a primitive root of unity of degree $p$.

Then by Theorem~\ref{thm:Verlinde} and
Proposition~\ref{Prop:defarg}
we have 
\begin{multline}\label{eq:part-of-inequ}
    \dim\rep{p}{r}[W(\vect N^{i^1,j^1},\dots,\vect
  N^{i^k,j^k};0)]\geq
  \dim\rep{p}{r}[W(\vect N^{i^1,j^1},\dots,\vect
  N^{i^k,j^k};\va)]\geq
  d_{p,r}[i^1,j^1;\dots;i^k,j^k]
\end{multline}
for any $\va \in \oC^k$.

One aim of this paper is to show that (\ref{eq:part-of-inequ}) is indeed
an equality 
for arbitrary vector~$\va \in \oC^k$ (see Theorem~\ref{thm:Verlinde-1}).

\bigskip

It is easy to see from~\eqref{eq:def-d} that 
$d_{p,r}[i^1,j^1;i^2,j^2;\dots;i^k,j^k]$ is a sum with the step~$p$
of supernomial coefficients.
To write this sum
introduce the $m\times m$
matrix\,\footnote{See~\cite{Sch-Warn}}
\begin{equation}
  T_m=\tiny
  \begin{pmatrix}
   2&-1& 0& 0&\dots& 0& 0& 0& 0\\
  -1& 2&-1& 0&\dots& 0& 0& 0& 0\\
   0&-1& 2&-1&\dots& 0& 0& 0& 0\\
\dots&\dots&\dots&\dots&\dots&\dots&\dots&\dots&\dots\\
   0& 0& 0& 0&\dots&-1& 2&-1& 0\\
   0& 0& 0& 0&\dots& 0&-1& 2&-1\\
   0& 0& 0& 0&\dots& 0& 0&-1& 1
  \end{pmatrix}
\end{equation}
so $\left(T^{-1}_m\right)_{ij}=\min\{i,j\}$ for $1\leq i,j\leq m$.

With the vector $\vect{N}=(N_+,N_-;N_0,N_1,\dots,N_{p-2})$ we associate 
the other vector 
$$\vect{N'}=(N_0,N_1,\dots,N_{p-2},N_++N_-).$$
Define the vector~$\vect L=(L_1,L_2,\dots,L_{p})$ as follows
\begin{equation}\label{eq:def-L}
  \vect{L}=\vect{N'}T_{p}\,.
\end{equation}
\begin{Prop} Let $\vect N=\sum_{s=1}^k\vect N^{i^s,j^s}$ and $\vect L = 
\vect{N'}T_{p}$. Then
\begin{equation}\label{eq:supernom-sum-d}
 d_{p,r}[i^1,j^1;i^2,j^2;\dots;i^k,j^k]=
  \sum_{a\in \oZ} \qsup{\LL}{pa+N_--r}{}^{p+1}\,,
\end{equation}
\end{Prop}

In the case where at least one of~$i^s=p-1$
sum~\eqref{eq:supernom-sum-d} can be easily calculated. 
\begin{Rem} We have
  \begin{equation}
    d_{p,r}[i^1,j^1;i^2,j^2;\dots;i^k,j^k]=2^{L_1}\,3^{L_2}\,
            \dots\,(p-1)^{L_{p-2}} p^{L_{p-1}-1} (p+1)^{L_p}
  \end{equation}
whenever $L_{p-1}\geq1$ or equivalently at least 
one~$i^s$ is equal to~$p-1$.
\end{Rem}

At last let us describe the class of subalgebras $W(\vect N)$ which can 
be deformed in this way.

\begin{Prop}
Suppose that all the components of the vector
$\vect{L}=\vect{N'}T_{p-1}$ are non-negative.
Then we have
$$W(\vect N)\simeq W(\vect N^{i^1, j^1},\vect N^{i^2, j^2},\dots,\vect
N^{i^k, j^k};0)$$
for some $i^1, \dots, i^k;\, j^1, \dots, j^k$.
\end{Prop}

\section{Degenerations of $\aB_p$ and the main 
inequality\label{sec:algebras}}
\subsection{}
Let $I$ be a set with $\ell$ elements.
Let~$\vect{e}_a$ with~$a\in I$ be the standard basis in~$\oZ^\ell$.
We write elements~$\vect{x}$ from~$\oZ^\ell$ as rows.
Let~$\vect{x}\cdot\vect{y}$ denote the standard scalar product
in~$\oZ^\ell$. 

Let~$\Gamma(A)$ be the lattice with the
basis~$\vect{e}_a$ and the scalar product determined by
a symmetric matrix~$A$
\begin{equation}
  \scalar{\vect e_a}{\vect e_b}=A_{ab}\in\oZ\,,\qquad a,b\in I\,.
\end{equation}

Consider the algebra~$\bB(A)$ generated by vertex 
operators~$V(\vect e_a,z)$, $a\in I$. Emphasize that we do not involve 
vertex operators corresponding to $-\vect e_a$. 
Let us fix the decomposition of~$V(\vect e_a,z)$ into
formal power series
\begin{equation}\label{eq:vertex-modes-decomposition}
  V(\vect e_a,z)=\sum_{n\in\oZ}V(\vect e_a)_n\,z^{-n-1}\,.
\end{equation}

The algebra~$\bB(A)$ is naturally $\Gamma(A)$-graded.
Decomposition~\eqref{eq:vertex-modes-decomposition}
determines another one grading by operator~$\cL_0$
\begin{equation}
  \commut{\cL_0}{V(\vect e_a)_n}=-n V(\vect e_a)_n\,.
\end{equation}

This algebra acts in the space
\begin{equation}
  \mathsf{H}=\bigoplus_{\vect y \in\Gamma^\vee(A)}
    \mathsf{H}({\vect y})\,,
\end{equation}
where $\mathsf{H}({\vect y})$ is the Fock module with the
momentum~$\vect y$ and~$\Gamma^\vee(A)$ is the lattice dual 
to~$\Gamma(A)$.
The representation~$\mathsf{H}$ 
contains the subrepresentation~$\Brep{p}{\vect y}$
generated from the highest weight vector~$\ket{\vect y}$ 
in~$\mathsf{H}({\vect y})$. 
The representation $\Brep{p}{\vect y}$  plays a role of Verma module,
namely, it is maximal among representation with the same annihilation
conditions satisfied by the vector~$\ket{\vect y}$. 
Clearly $\Brep{p}{\vect y}$ admits the $\Gamma^\vee(A)$-grading; to fix it
let us set~$\deg\ket{\vect y}=0$.

The dual space~$\Brep{p}{\vect y}^*$ can be identified with the space of
correlation functions. Namely, any homogeneous 
linear functional~$\bra{\theta}\in\Brep{p}{\vect y}^*$ of degree $\sum
n_i\vect e_i$ determines and is determined by the function
\begin{multline}
  Y_{\theta,\vect y}(x^1_1,x^1_2,\dots ,x^1_{n_1};
  x^2_1,x^2_2,\dots, x^2_{n_2};\dots;
  x^\ell_1,x^\ell_2,\dots x^\ell_{n_\ell})=\\
  \bra{\theta}V(\vect e_1,x^1_1)V(\vect e_1,x^1_2)\dots
       V(\vect e_1,x^1_{n_1})V(\vect e_2,x^2_1)V(\vect e_2,x^2_2)\dots
       V(\vect e_2,x^2_{n_2})\dots\\
        V(\vect e_\ell,x^\ell_1)V(\vect e_\ell,x^\ell_2)\dots
       V(\vect e_\ell,x^\ell_{n_\ell}) \ket{\vect y}.
\end{multline}

\begin{Prop}\label{prop:corr-funcs}  
Correlation functions $Y_{\theta,\vect y}$ are products of the polynomial
\begin{equation}
  (x^1_1x^1_2\dots x^1_{n_1})^{w_1}
  (x^2_1x^2_2\dots x^2_{n_2})^{w_2}\dots
  (x^\ell_1x^\ell_2\dots x^\ell_{n_\ell})^{w_\ell}
        \prod\limits_{\atop{1\leq a\leq\ell}{1\leq i<j\leq n_a}}
    (x^a_i-x^a_j)^{A_{aa}}
  \prod\limits_{\atop{1\leq a<b\leq\ell}{\atop{1\leq i\leq n_a}{1\leq
        j\leq n_b}}} (x^a_i-x^b_j)^{A_{ab}}
\end{equation}
and a polynomial symmetric in each group of 
variables~$x^a_1,x^a_2,\dots,x^a_{n_a}$. 
The numbers~$w_i$ are determined by the annihilation conditions
satisfied by the vector~$\vect y$.
\end{Prop}

In other words, we have the decomposition
\begin{equation}
  \Brep{p}{\vect y}^*=\bigoplus_{\vect x\in\Gamma(A)}\Brep{p}{\vect
  y}^*(\vect x)
\end{equation}
and $\Brep{p}{\vect y}^*(\vect x)$ can be identified with the space of
correlation functions of degree $\vect x$.
Using Proposition~\ref{prop:corr-funcs} we obtain the following
formula of Gordon type. 
\begin{Prop}
The character of the space $\Brep{p}{\vect y}$ is given by the formula
\begin{equation}
  \chrct\Brep{p}{\vect y}(q,\vect z)=
   \sum_{\vect{n}\in\oZ^\ell}
  \vect{z}^{\vect{n}}\,
  q^{\half\vect{n}A\cdot\vect{n} + \vect{v}\cdot\vect{n}}
 \prod\limits_{a\in I}
 \frac{1}{(q)_{\vect{e}_a\cdot\vect{n}}}
\end{equation}
where $\vect{v}=-\half\diag(A)+(w_1,w_2,\dots,w_\ell)$
and~$\vect{z}^{\vect{n}}=
z_1^{n_1}z_2^{n_2}\cdots z_\ell^{n_\ell}$.
\end{Prop}

Let $\vect N$ be an $\ell$-dimensional
vector with integer components.
Define the subspaces $W^A(\vect N) \subset \bB(A)$ as follows:
\begin{equation}
  W^A(\vect N) \ \mbox{is spanned by}\ 
\{V(\vect e_a)_n \ |\   n\geq N_a,\ a\in I\}.
\end{equation}
Then $W^A(\vect N)$-coinvariants of~$\Brep{p}{\vect y}^*$ can be easily
calculated in terms of correlation functions.  Namely the
dual space consists of functions given in 
Proposition~\ref{prop:corr-funcs}
with a restriction on the degrees in each variable.
This gives the following 
\begin{Prop}\label{Prop:char} The character 
of~$\Brep{p}{\vect y}[W^A(\vect N)]=\Brep{p}{\vect y}^*/
                    W^A(\vect N)\Brep{p}{\vect y}^*$
is given by the formula
\begin{equation}
  \chrct\Brep{p}{\vect y}[W^A(\vect N)](q,\vect z)=
   \sum_{\vect{n}\in\oZ^\ell}
  \vect{z}^{\vect{n}}\,
  q^{\half\vect{n}A\cdot\vect{n} + \vect{v}\cdot\vect{n}}
 \prod\limits_{a\in I}
 \qbin{\vect{e}_a\cdot(\vect{N}+\vect{n}-\vect{n}A+\vect{w})}
                   {\vect{e}_a\cdot\vect{n}}{q}
\end{equation}
where $\vect{w}=(w_1,w_2,\dots,w_\ell)$.
\end{Prop}

\subsection{}
Consider the family of algebras $\bA^d_p$ obtained from the 
algebra~$\aB_p$ by the following recursive procedure.

Introduce a filtration on~$\aB_p$ by assigning degree 1 to the
generators~$\genp{p}$  and~$\genm{p}$.
The associated graded algebra~$\bA^{0}_p$ is generated 
by~$\bgenp{p}$, $\bgenm{p}$. Define the currents~$\bgenh{p}{m}$
by the following OPE:
$$\bgenp{p}(z)\bgenm{p}(w)=\frac{\bgenh{p}{0}}{(z-w)^{p-1}}+\,\dots\,+
  \frac{\bgenh{p}{p-2}}{z-w}+\mbox{\rm regular terms}.$$
Note that the original algebra~$\aB_p$ is the central extension 
of~$\bA_p^0$.

At the next step we assign degrees~$1$ to~$\bgenp{p}$, $\bgenm{p}$
and~$\bgenh{p}{0}$. 
After taking graded object we obtain the 
algebra~$\bA^{1}_p$ whose generators we denote abusing notations again 
by~$\bgenp{p}$, $\bgenm{p}$ and~$\bgenh{p}{0}$.
The algebra~$\bA^{1}_p$ is generated by
currents~$\bgenp{p}$, $\bgenm{p}$ and~$\bgenh{p}{0}$ that have OPEs
with the following properties
\begin{align}
  \bgenpm{p}(z)\bgenpm{p}(w)&\sim(z-w)^p\,,\qquad
  \bgenpm{p}(z)\bgenh{p}{0}(w)\sim(z-w)^1\,,\qquad
  \bgenh{p}{0}(z)\bgenh{p}{0}(z)\sim(z-w)^2\\
  \bgenp{p}(z)\bgenm{p}(w)&=\frac{\bgenh{p}{1}}{(z-w)^{p-2}}+\,\dots\,+
  \frac{\bgenh{p}{p-2}}{z-w}+\mbox{\rm regular terms}
\end{align}

To make the step~$\bA^{d-1}_p\to\bA^{d}_p$ we introduce a filtration by
attaching degree~$1$ to~$\bgenp{p}$, $\bgenm{p}$ 
and~$\bgenh{p}{0}$, \dots, $\bgenh{p}{d-1}$.
As the result we obtain the algebra~$\bA^{d}_p$
with generators~$\bgenp{p}$, $\bgenm{p}$ and~$\bgenh{p}{i}$, $0\leq
i\leq d-1$ that have the following OPEs
\begin{align}\label{eq:OPEs-at-d}
\bgenpm{p}(z)\bgenpm{p}(w)&\sim(z-w)^p\,,\qquad
  \bgenpm{p}(z)\bgenh{p}{m}(w)\sim(z-w)^{\theA_{\pm m}}\,,\kern10pt
\bgenh{p}{n}(z)\bgenh{p}{m}(w)\sim(z-w)^{\theA_{nm}}\\
  \bgenp{p}(z)\bgenm{p}(w)&=\frac{\bgenh{p}{d}}{(z-w)^{p-d-1}}+\,\dots\,+
  \frac{\bgenh{p}{p-2}}{z-w}+\mbox{\rm regular terms}
\end{align}
determined by the 
following~$(d+2)\times(d+2)$ matrix
\begin{equation}\label{eq:the-matrix}
 \theA=\tiny
\begin{pmatrix}
     p &-p+d+1 &  1  &  2  &  3  &\dots&  d\\
     -p+d+1 & p &  1  &  2  &  3  &\dots&  d\\
       1  &  1  &  2  &  2  &  2  &\dots&   2 \\
       2  &  2  &  2  &  4  &  4  &\dots&   4 \\
       3  &  3  &  2  &  4  &  6  &\dots&   6 \\
     \dots&\dots&\dots&\dots&\dots&\dots&\dots\\
      d & d &  2  &  4  &  6  &\dots& 2d
  \end{pmatrix}
\end{equation}
where the first row and column corresponds to the index $+$ 
the second ones to $-$ the third to $0$ and so on till $d-1$.
\begin{Prop}
There is the surjective homomorphism 
\begin{equation}\label{eq:map}
  \bB(\theA)\to\bA^d_p
\end{equation}
acting on generators as
\begin{equation}\label{eq:map-currents}
  V(\vect e_+,z)\to\bgenp{p}(z)\,,\qquad
  V(\vect e_-,z)\to\bgenm{p}(z)\,,\qquad
  V(\vect e_i,z)\to\bgenh{p}{i}(z)\,,\quad0\leq i\leq d-1\,.
\end{equation}
\end{Prop}
\begin{proof}
From~\eqref{eq:OPEs-at-d} we can deduce that all relations
between~$V(\vect e_a,z)$ hold for their images in~$\bA^d_p$. 
So~\eqref{eq:map} defines a homomorphism. As the image
contains all the generators of $\bA^d_p$, this map is surjective.
\end{proof}

Note that the map~\eqref{eq:map} is graded. Namely, it
preserves the grading with respect to $\cL_0$ and specialize the
$\Gamma$-grading of $\bB(\theA)$ into $\genh{p}{0}$-grading of $\bA^d_p$. 
This specialization is defined by the functional $\vect u$ on $\Gamma$
such that $\vect u(\vect e_+) = 1$, $\vect u(\vect e_-) = -1$, $\vect
u (\vect e_i) =0$ (see~\eqref{eq:vect-u}).

Any representation $\rep{p}{r}$ of~$\aB_p$ contains the highest weight
vector. Therefore the filtration in~$\aB_p$ determines a filtration
in~$\rep{p}{r}$ and associated graded algebra~$\bA_p^0$
acts in the associated graded module~$\brep{p}{r}^0$. We can 
repeat this procedure with $\brep{p}{r}^0$ and the image of the
highest weight vector.
In this way we obtain the family of representations~$\brep{p}{r}^d$
of algebras~$\bA_p^d$.
Evidently characters of all~$\brep{p}{r}^d$, $0\leq d\leq p-2$
are equal to each other and equal to character
of~$\rep{p}{r}$ given by ~\eqref{eq:char-I}.

Note that~\eqref{eq:map} allows us to treat any $\bA^d_p$-module
as a $\bB(\theA)$-module. According to this we have
\begin{Prop} Let~$\vect y=r\vect u$. Then 
there is the surjective mapping of $\bB(\theA)$-modules
\begin{equation}\label{eq:map-reps} 
\Brep{p}{\vect y}\to\brep{p}{r}^d
\end{equation}
mapping the vector $\ket{\vect y}$ into the highest
weight vector.
\end{Prop} 
\begin{proof}
It follows from the observation that the highest
weight vector of~$ \brep{p}{r}^d$ satisfy the same annihilation
conditions as~$\ket{\vect y}$.
\end{proof}

We will show (see Corollary~\ref{Cor:iso}) that~\eqref{eq:map-reps} is
indeed an isomorphism. Note that this map is bi-graded with respect to
$\cL_0$ and $\vect u (\Gamma)$.

Define the subspaces $W^d (\vect N) \subset \bA^d_p$ recursively 
from $W(\vect N)$
as the adjoint graded subspaces. 

\begin{Prop}\label{Prop:map-sub}\
Let $0 \le d \le p-1$. Then 
$W^d (\vect N)$ contains
the image of $W^{\theA}(\vect N)$ under the map~\eqref{eq:map}.
\end{Prop}
\begin{proof}
  Note that $W(\vect N)$ contains linear combinations of the form
$$
\genh{p}{i}_s+c_1\genh{p}{i-1}_s+\dots+c_i\genh{p}{0}_s+c_{i+1}
$$
for $s\geq N_i-i$. Therefore $W^d(\vect N)$ contains $\bgenh{p}{i}_s$
for $s\geq N_i-i$ hence (according to~\eqref{eq:modes-h} 
and~\eqref{eq:vertex-modes-decomposition})
it contains the image of~$W^{\theA}(\vect N)$.
\end{proof}

The most important case of this statement is $d=p-1$, where the
condition on $\vect N$ is empty.

\begin{Cor} Under the conditions of Proposition~\ref{Prop:map-sub}
we have the surjective map
\begin{equation}
  \label{eq:sur-coinv}
  \Brep{p}{r\vect u}[W^{\theA}(\vect N)]\to\brep{p}{r}^d[W(\vect N)]
\end{equation}
\end{Cor}

This map is an isomorphism only under some additional conditions (see
Corollary~\ref{Cor:iso}).

\def \gr {{\rm gr}\,}

The deformation argument (see for example~\cite{FKLMM2}) implies that 
for any filtered representation $V$ of filtered (sub)algebra $W$ we
have 
$$\dim (\gr V [ \gr W]) \ge \dim V[W].$$

In our case it implies
\begin{equation}
  \label{eq:half-of-inequality}
  \dim\brep{p}{r}^d[W^d(\vect N)]\geq\dim\brep{p}{r}^{d-1}[W^{d-1}(\vect N)]
  \geq\dots\geq\dim\brep{p}{r}^0[W^0(\vect N)]\geq\dim\rep{p}{r}[W(\vect N)]
\end{equation}
Taking \eqref{eq:sur-coinv}, \eqref{eq:half-of-inequality}
and Proposition~\ref{eq:part-of-inequ} together
we obtain the {\sl main inequality}
\begin{Prop} For $0 \le d \le p-1$ we have
\begin{multline}
  \label{eq:the-main-inequality}
  \dim\Brep{p}{r\vect u}[W^{\theA}(\vect N)]\geq
  \dim\rep{p}{r}[W(\vect N^{i^1,j^1},\vect N^{i^2,j^2},\dots,\vect
  N^{i^k,j^k};0)]\geq\\
  \dim\rep{p}{r}[W(\vect N^{i^1,j^1},\vect N^{i^2,j^2},\dots,\vect
  N^{i^k,j^k};\vect a)]\geq
    d_{p,r}[i^1,j^1;i^2,j^2;\dots;i^k,j^k]
  \end{multline}
\end{Prop}

\section{Combinatorics of the characters\label{sec:combinator}}
\subsection{}
Following \cite{Sch-Warn} introduce $q$-supernomial coefficients as
follows
\begin{equation}\label{sl2-supernomial}
\qsup{\LL}{a}{q}^{k+1} =
\sum\limits_{\atop{n_1,\dots, n_k  \in \oZ}{n_1 + \dots + n_k =a}}\, 
{q}^{\sum\limits_{i=2}^{k}
n_{i-1}(-n_i+\sum\limits_{j=i}^{k}L_j)}
    \qbin{L_k}{n_k}{q}\qbin{L_{k-1}+n_{k}}{n_{k-1}}{q}
     \dots
  \qbin{L_{1}+n_{2}}{n_{1}}{q},
\end{equation}
where $\LL=(L_1,\,L_2,\,\dots,\,L_{k})$ is a $k$-dimensional vector
with nonnegative integer entries. 
The sum in the RHS is a well-defined polynomial in~$q$
because $q$-binomial coefficients
vanish whenever the bottom index becomes greater than the top one.

Let $\NN = (N_+, N_-; N_0, \dots, N_{d-1})$ be a vector in $\oZ^{d+2}$; 
let $\NN' = (N_0, \dots, N_{d-1}, N_+ +
N_-)$ and let $\LL = \NN' T_{d+1}^{-1}$.
The aim of this section is to prove the following
polynomial identity.
\begin{Thm}\label{ta} Let $0 \le d < 2p-2$.
Suppose that $L_i \ge 0$ and additionally
\begin{equation}\label{wellb}
2N_+ - N_{d-1} \geq -(2p-d-2), \qquad 2N_- - N_{d-1} \geq -(2p-d-2). 
\end{equation}
Then
\begin{equation}\label{cmain}
\sum_{\vect{n}\in\oZ^\ell}
  z^{n_+ - n_-}\,
  q^{\half\vect{n}{\theA}\cdot\vect{n}}
 \prod\limits_{a\in I}
 \qbin{\vect{e}_a\cdot(\vect{N}+\vect{n}-\vect{n}{\theA})}
                   {\vect{e}_a\cdot\vect{n}}{q}
= \sum_{a\in \oZ} z^a q^{pa^2/2} \qsup{\LL}{pa+N_-}{q}^{d+2}.
\end{equation}
where $I=\{+,-;0,1,\dots,d-1\}$
\end{Thm}

We call (\ref{wellb}) well-balanced conditions. Without them identity
(\ref{cmain}) can fails.

To eliminate the well-balanced conditions we use an improved version of
$q$-binomial coefficients:
\begin{equation}\label{eq:qbinplus}
\qbin{n}{m}{q}^+ = \left\{
\begin{array}{ll}
\qbin{n}{m}{q} & n\ge 0\\
(-1)^{n-m} q^{-\frac{(n-m)^2+ (n-m)}2} \qbin{-m-1}{-n-1}{q^{-1}} & n<0
\end{array}\right.
\end{equation}
In particular, for $q\to 1$ the integer $\qbin{n}{m}{1}^+$ is the 
coefficient before~$z^{-m}$ in the Laurent expansion 
of~$(1+z^{-1})^n$ at~$z=0$.

This definition is motivated by the following relations.

\begin{Prop}\label{Prop:exqbin}
The polynomials $\qbin{n}{m}{q}^+$ are uniquely defined by the
$q$-Pascal identities
\begin{align}\label{qpascal2}
\qbin{n}{m}{q}^+ & = \  q^m \qbin{n-1}{m}{q}^+ +  \qbin{n-1}{m-1}{q}^+
& \mbox{for any \ $n$, $m$}\\
\qbin{n}{m}{q}^+ & = \qbin{n-1}{m}{q}^+ + q^{n-m} \qbin{n-1}{m-1}{q}^+ 
& \mbox{for any \ $n$, $m$}
\label{qpascal}\end{align} 
and conditions
\begin{equation}\label{eq:bound}
\qbin{m}{m}{q}^+ = 1, \qquad \qbin{n}{0}{q}^+ = 0\ \ \mbox{for $n<0$}
\end{equation}
which need to  be checked only for certain $m$ and $n$.
\end{Prop}
\begin{proof}
Indeed any set of polynomials $\qbin{n}{m}{q}^?$ satisfying
\eqref{qpascal2} and \eqref{qpascal} can be expressed in the following
way:
$$\qbin{n}{m}{q}^? = \left\{
\begin{array}{ll}
P_1 \cdot (-1)^{n-m} q^{-\frac{(n-m)^2+ (n-m)}2}
\qbin{-m-1}{n-m}{q^{-1}} &  n<0,\ m<0\\
P_2 \cdot (-1)^m q^{-\frac{m^2 + m}2} \qbin{-(n-m)-1}{m}{q^{-1}} & n<0,\ m\ge
0\\
(P_1 + P_2) \cdot \qbin{n}{m}{q} & n\ge 0
\end{array}\right.\quad
\begin{picture}(80,30)(-10,40)
\epsfig{file=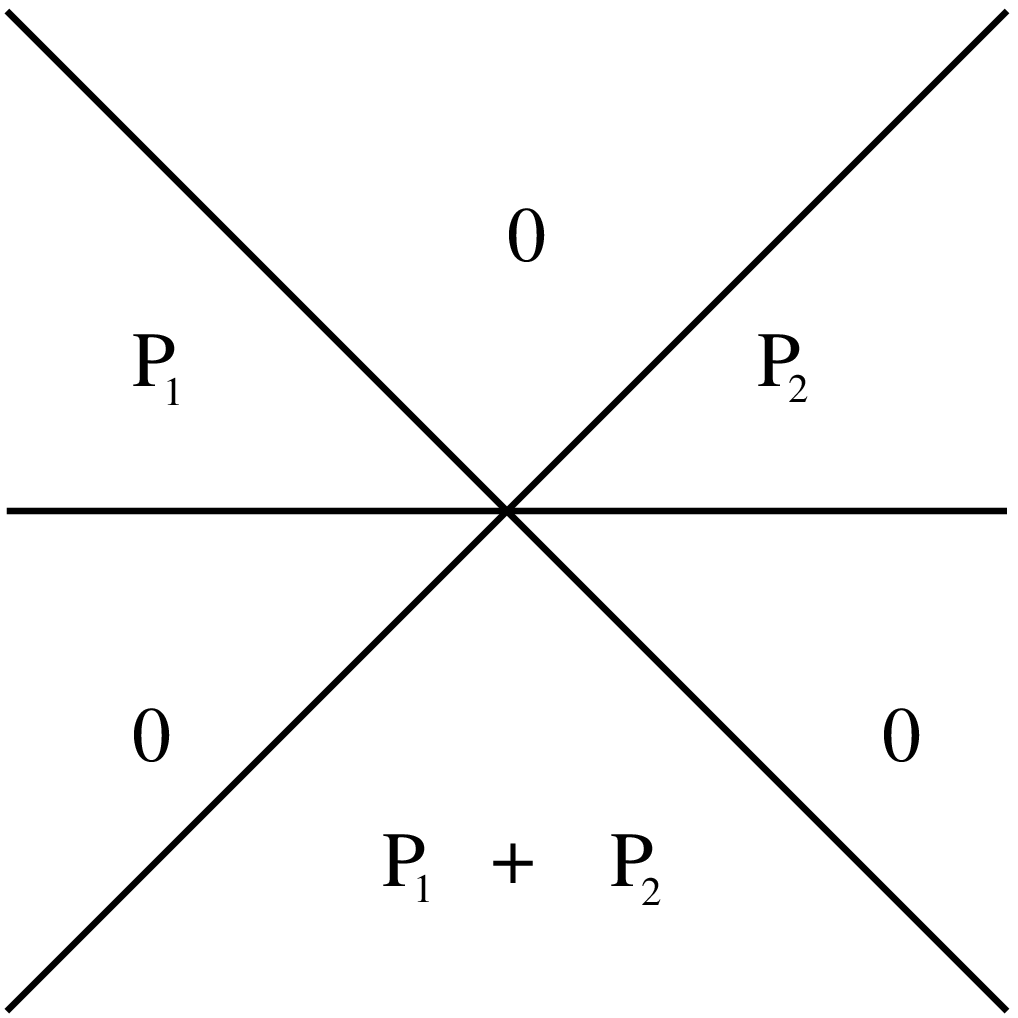, height=80pt}
\end{picture}
$$
where $P_1$ and $P_2$ are certain polynomials in $q$.

In our case we need only to check that $P_1 =1$ and $P_2 =0$.
\end{proof}

Note that for usual $q$-binomial coefficients \eqref{qpascal2} and
\eqref{qpascal} fails for $n=m=0$.

\begin{Thm}\label{tb} Let $0 \le d < 2p-2$.
Suppose that $L_i \ge 0$. Then 
\begin{equation}\label{cimp}
\sum_{\vect{n}\in\oZ^\ell}
  z^{n_+ - n_-}\,  
  q^{\half\vect{n}{\theA}\cdot\vect{n}}                        
 \prod\limits_{a\in I}
 \qbin{\vect{e}_a\cdot(\vect{N}+\vect{n}-\vect{n}{\theA})}
                   {\vect{e}_a\cdot\vect{n}}{q}^+
= \sum_{a\in \oZ} z^a q^{pa^2/2} \qsup{\LL}{pa+N_-}{q}^{d+2}.
\end{equation}
\end{Thm}
Emphasize that in the left hand side the sum goes over all integer vectors.

\subsection{}
Here we show that the left hand side of (\ref{cimp}) is a well
defined polynomial
and that theorem~\ref{tb} implies theorem~\ref{ta}.

\def \nn {\vect n}

Let $P_{\NN}(\nn)=  \prod\limits_{a\in I}
 \qbin{\vect{e}_a\cdot(\vect{N}+\vect{n}-\vect{n}{\theA})}
                   {\vect{e}_a\cdot\vect{n}}{q}^+$. We need to prove that
if $L_i \ge 0$ for any $i$ then there is only finite number of vectors
$\nn$ such that $P_{\NN}(\nn) \ne 0$. And if additionally we impose the
well-balance conditions (\ref{wellb}) then $P_{\NN}(\nn) \ne 0$ implies
$n_a \ge 0$, $a = +, -, 0, 1, \dots, d-1$, 
where~$n_a=\vect e_a\cdot\vect n$. 

\begin{Prop}
If $P_{\NN}(\nn) \ne 0$ then
\begin{align}\label{rall}
N_a\ge\vect e_a\cdot\nn {\theA} &\qquad 
   \mbox{for any}\ a = +,-,0,1,2,\dots, d-1;\\
\label{rneg}
N_a < \vect e_a\cdot(\nn {\theA}-\nn) &\qquad \mbox{if}\ n_a < 0.
\end{align}
\end{Prop}
\begin{proof} Condition~\eqref{rall} follows from the fact that
$\qbin{n}{m}{q}^+ = 0$ if $n<m$. Condition~\eqref{rneg} follows from the
fact that 
$\qbin{n}{m}{q}^+ = 0$ if $m<0$ and $n\ge 0$.
\end{proof}

\begin{Lemma}
Let $1\le i \le d$.
If $L_i \ge 0$ and $n_{i-1} <0$ then $P_{\NN}(\nn) = 0$.
\end{Lemma}
\begin{proof} Suppose that $n_{i-1} <0$ and $P_{\NN}(\nn)
\ne 0$.  
Then
adding~\eqref{rneg} for $a=i-1$ twice and subtracting \eqref{rall} for
$a=i-2$ (if $i>1$) and $a=i$ (if $i<d$; if $i=d$ then for $a=+$ and
$a=-$) we obtain $L_i<0$. 
\end{proof}

\begin{Lemma}
If $L_{d+1} \ge 0$ and $n_+ <0$, $n_-<0$ then $P_{\NN}(\nn) = 0$.
\end{Lemma}
\begin{proof}
Proof is similar. Here we add~\eqref{rneg} for $a=+$ and for $a=-$
then subtract \eqref{rall} for
$a=d-1$. 
\end{proof}

\begin{Lemma}
Suppose that $n_+<0$ and $n_-\ge 0$. Then $P_{\NN}(\nn) = 0$ for 
$$n_+
\le C_+=(2N_+ - N_{d-1}) / (2p-d-2).$$
\end{Lemma}
\begin{proof}
Suppose that $n_+<0$ and $P_{\NN}(\nn) \ne 0$. 
Then
adding~\eqref{rneg} for $a=+$ twice and subtracting~\eqref{rall} for
$a=d-1$
we obtain
$$ 2N_+ - N_{d-1} < (2p-d-2) n_+ - (2p-d-2) n_-.$$
Taking into account $n_-\ge 0$ and $2p-d-2>0$ we deduce $(2p-d-2) n_+ >
2N_+ -
N_{d-1}$.
\end{proof}

Clearly, we have the similar statement for $n_- < 0$ and $n_+ \ge 0$.

So if $L_i \ge 0$ for any $i$ and the well-balanced conditions are
satisfied then $P_\NN(\nn)\ne 0$ only if  $n_a \ge 0$ for any $a =
+,-,0,\dots,d-1$, hence
theorem~\ref{tb} implies theorem~\ref{ta}.

Under the conditions of theorem~\ref{tb} we have $n_i \ge 0$ for $i=0
\dots d-1$ and $n_+> C_+$,  $n_- > C_-$ otherwise~$P_{\NN}(\nn) = 0$.
Adding (\ref{rall}) for $a = +$ and $a=-$, we obtain
$$N_+ + N_- \ge (d+1) n_+ + (d+1) n_- + 2n_0 + 4 n_1 + 6 n_2 + \dots.$$
As all the coefficient in the right hand side of this inequality are
positive, there exists only finite number of such vectors $\nn$.

\subsection{}
Here we
compare the left hand side $\chi^L_{p,d}[\NN](q,z)$ and the
right hand side $\chi^R_{p,d}[\LL,N_-](q,z)$ of (\ref{cimp}).

Note that for $q$-supernomial coefficients there is an identity
$$\qsup{\LL+\ve_k}{a}{q}^{k+1} = \qsup{\LL+\ve_{k-1}}{a-1}{q}^{k+1} +
q^a \qsup{\LL}{a}{q}^{k}.$$

Therefore, for the right hand side we have
\begin{equation}\label{rec1}
\chi^R_{p,d}[\LL+ \ve_{d+1},N_-]=\chi^R_{p,d}[\LL+
\ve_{d},N_- -1] + z^{-1} \,q^{N_- -p/2}\, \chi^R_{p,d}[\LL, N_- -
p].
\end{equation}

Also as $\qsup{(L_1, \dots, L_{k-1},0)}{a}{q}^{k+1} = \qsup{(L_1, \dots,
L_{k-1})}{a}{q}^{k}$, we have 
\begin{equation}\label{rec2}
\chi^R_{p,d}[(L_1, \dots, L_d,0), N_-] = \chi^R_{p,d-1}[(L_1, \dots,
L_d), N_-].
\end{equation}

The idea of the proof is to check (\ref{rec1}) and (\ref{rec2}) for the
left hand side.
\begin{Prop}
Let $A$ be a symmetric $m\times m$ matrix with integer components, let
$\vu$, $\NN$ be vectors in $\oZ^m$.
Consider the polynomial
\begin{equation}
  \chi_{A,\vect{u}}[\vect{N}](q,z)=
   \sum_{\vect{n}\in\oZ^m}
  z^{\vect{u}\cdot\vect{n}}\,
  q^{\half\vect{n}A\cdot\vect{n}}
 \prod\limits_{a}
 \qbin{\vect{e}_a\cdot(\vect{N}+\vect{n}-\vect{n}A)}
                   {\vect{e}_a\cdot\vect{n}}{q}^+.
\end{equation}
Then for any $a$ we have
\begin{equation}\label{eq:main-recurrence}
  \chi_{A,\vect{u}}[\vect{N}]=
  \chi_{A,\vect{u}}[\vect{N}-\vect{e}_a]
 +z^{u_a}\,
  q^{N_a - A_{aa}/2}
\, \chi_{A,\vect{u},\vect{v},\vect{w}}[\vect{N}-\vect{e}_aA].
\end{equation}
\end{Prop}
\begin{proof}
Identity~\eqref{eq:main-recurrence} follows from the
identity~\eqref{qpascal} applied to the factor
$\qbin{\vect{e}_a\cdot(\vect{N}+\vect{n}-\vect{n}A)}
                   {\vect{e}_a\cdot\vect{n}}{q}^+$ in each summand.
\end{proof}

Applying (\ref{eq:main-recurrence}) with $a=-$ to
$\chi^L_{p,d}[\NN]$, we check (\ref{rec1}) for the left hand side.
Note that for $a=0,\dots, d-1$ (\ref{eq:main-recurrence}) implies recurrence
relations for the $q$-supernomial coefficients same as in \cite{Sch-Warn}.

To check (\ref{rec1}) for the left hand side we need the following
identity.
\begin{Lemma}
There is the following relation
\begin{equation}\label{rdc}
  \qbin{N}{n}{q}^+
  \qbin{M}{m}{q}^+=
  \sum_{l\in\oZ}q^{(n-l)(m-l)}
        \qbin{N - m}{n - l}{q}^+
        \qbin{M - n}{m - l}{q}^+
        \qbin{N + M - n - m + l}{l}{q}^+
\end{equation}
\end{Lemma}
\begin{proof}
First let us check that the sum in the right hand side is indeed
finite.

Note that if $N+M-n-m<0$ then the third factor 
is always zero, so entire sum is zero. If $N+M-n-m\ge 0$ then either
$N-m \ge 0$ or $M-n \ge 0$. In the first case we have non-zero
summands only if $n+m-N \le l \le n$, in the second case only if
$n+m-M \le l \le m$.

\def \rh {{\rm RHS}}

Let us denote the polynomial in the right hand side by $\rh^{N,M}_{n,m}$.
The idea of proof is to apply Proposition~\ref{Prop:exqbin} to
$\rh^{N,M}_{n,m}$ as a set of polynomials indexed by $N$ and $n$. It
means that we prove that
\begin{equation}\label{eq:biqp}
\rh^{N,M}_{n,m} = q^n \rh^{N-1,M}_{n,m} + \rh^{N-1,M}_{n-1,m},\quad
\rh^{N,M}_{n,m} = \rh^{N-1,M}_{n,m} + q^{N-n} \rh^{N-1,M}_{n-1,m}
\end{equation}
for any $N$,$n$,$M$,$m$ and
\begin{equation}\label{eq:bibnd}
\rh^{n,M}_{n,m} = \qbin{M}{m}{q}^+ \qquad \rh^{N,M}_{0,m} = 0
\end{equation}
for any $M$, $m$ and certain $N<0$, $n$.

\newcommand\sft[6]{\left[
\begin{array}{rrr}
#1 & #2 & #3 \\
#4 & #5 & #6
\end{array}\right]}

To prove \eqref{eq:biqp} fix for a moment $N$, $M$, $n$, $m$ and
introduce the notation 
$$\sft{a}{b}{c}{d}{e}{f} = q^{(n-l)(m-l)}
        \qbin{N - m + a}{n - l + d}{q}^+
        \qbin{M - n + b}{m - l + e}{q}^+
        \qbin{N + M - n - m + l + c}{l + f}{q}^+.$$
Then substitution $l+1$ into $l$ leads to the identity
$$\sum_{l\in\oZ} \sft{a}{b}{c-1}{d+1}{e+1}{f-1} = \sum_{l\in\oZ} q^{2l-m-n+1}
\sft{a}{b}{c}{d}{e}{f}.$$
Using this identity together with \eqref{qpascal2} we obtain
\begin{multline*}
\rh^{N,M}_{n,m} = \sum \sft{0}{0}{0}{0}{0}{0} = \sum q^{n-l}
\sft{-1}{0}{0}{0}{0}{0} + \sum \sft{-1}{0}{0}{-1}{0}{0} =\\  
= \sum q^{n} \sft{-1}{0}{-1}{0}{0}{0} + \sum q^{n-l}
\sft{-1}{0}{-1}{0}{0}{-1} + \sum \sft{-1}{0}{0}{-1}{0}{0} = \\
= \sum q^{n} \sft{-1}{0}{-1}{0}{0}{0} + \sum q^{l-m}
\sft{-1}{0}{0}{-1}{-1}{0} + \sum \sft{-1}{0}{0}{-1}{0}{0} = \\
= \sum q^{n} \sft{-1}{0}{-1}{0}{0}{0} + \sum q^{l-m} \sft{-1}{1}{0}{-1}{0}{0} =
q^n \rh^{N-1,M}_{n,m} + \rh^{N-1,M}_{n-1,m},
\end{multline*}
where sums go over $l \in \oZ$.
Proof of the second part of \eqref{eq:biqp} is similar.

To check \eqref{eq:bibnd} we take $n=M$ and any negative $N<m-M$. We have only
one non-zero summand in $\rh^{M,M}_{M,m}$, namely  for $l=m$, and it
coincides with the $q$-binomial given in \eqref{eq:bibnd}. If $N<m-M$
and $n=0$ then $N+M-m-n <0$ and, as described above, the result is
zero.

So proposition~\eqref{Prop:exqbin} implies the lemma.
\end{proof}

Applying (\ref{rdc}) to the factors in $\chi^L_{p,d-1}[(N_+,N_-;N_0 \dots,
N_{d-2})]$ corresponding to $n_+$ and $n_-$ we obtain
$$\chi^L_{p,d-1}[(N_+,N_-;N_0 \dots,
N_{d-2})] = \chi^L_{p,d}[(N_+,N_-;N_0 \dots,
N_{d-2},N_++N_-)],$$
that is, (\ref{rec2}) for the left hand side.

Now we can prove theorem~\ref{tb} by induction on $d$ and $L_{d+1}$. To
complete the proof we need only to check the case $d=0$ manually.
\begin{Lemma} Suppose that $M+S \ge 0$. Then 
  \begin{equation}\label{eq:Knuth-id}
    \sum_{k\in\oZ}q^{k^2+ak}\qbin{M}{a+k}{q}^+\qbin{S}{k}{q}^+=
    \qbin{M+S}{S+a}{q}
  \end{equation}
\end{Lemma}
\begin{proof}
First let us check that the sum in the left hand side is indeed
finite.
As $M+S \ge 0$ we have either $M\ge 0$ or $S \ge 0$. In the first case
we have non-zero 
summands only if $-a \le k \le M-a$, in the second case only if
$0 \le k \le S$.

\def \lh {{\rm LHS}}

Let us denote the polynomial in the left hand side by $\lh^{M,S}_a$.
Clearly, we have $\lh^{M,S}_a = \lh^{S,M}_{-a}$ and the same identity
for the right hand side. Therefore we can restrict ourselves to the
case $S\ge 0$ and prove that 
\begin{equation}\label{eq:iknuth}
\lh^{M,S}_a =
\qbin{M+S}{S+a}{q}^+ \quad \mbox{for $S \ge 0$ and any $M$.}
\end{equation}

Applying \eqref{qpascal} to the second factor in $\lh^{M,S}_a$, we
obtain
$$\lh^{M,S}_a = \lh^{M,S-1}_a + q^{S+a}\lh^{M,S-1}_{a+1},$$
that is \eqref{qpascal2} for the right hand side of \eqref{eq:iknuth}.
Therefore it remains
to check \eqref{eq:iknuth} for $S=0$.

If $S=0$ then the only non-zero summand in $\lh^{M,S}_a$ is
for $k=0$ and it is equal to the right hand side of \eqref{eq:iknuth}
\end{proof}

The case $d=0$ of Theorem~\ref{tb} can be obtained by substituting
$M=N_+ -(p-1)(n_+-n_-)$,  $S=N_- +(p-1)(n_+-n_-)$, $a=n_+-n_-$
and~$k=n_-$ in~\eqref{eq:Knuth-id}.

\section{Dimension and character formulas for
the coinvariants\label{sec:main-result}}

\begin{Thm}\label{thm:Verlinde-1}\label{flatness}
For any $\va \in \oC^k$ we have
  \begin{equation}
    \dim \rep{p}{r}[W(\vect N^{i^1,j^1},\vect N^{i^2,j^2},\dots,\vect
N^{i^k,j^k};\va)]=d_{p,r}[i^1,j^1;i^2,j^2;\dots;i^k,j^k]\,.
  \end{equation}
\end{Thm}
\begin{proof}
First, note that due to the action of automorphisms $\spfn{\theta}$ we
can consider only coinvariants of the representation $\rep{p}{0}$.

The main idea of the proof is to check that the main
inequality~(\ref{eq:the-main-inequality}) is
indeed an equality. Here we have a chain of inequalities, where we know 
the first element by Proposition~\ref{Prop:char} and the last
element. 
As soon as we check that these two numbers are
equal we obtain that other numbers in (\ref{eq:the-main-inequality}) are
also equal.

Suppose that the vector $\NN$ satisfies the well-balance
conditions~(\ref{wellb}). Then this equality follows from theorem~\ref{ta}
by setting $z=1$ and $q=1$. But if conditions~(\ref{wellb}) are not
satisfied then this equality can be violated. 
So we should deduce general case
from the well-balanced one.

To this end consider the family of automorphisms  $\spfn{p\theta}$. They
preserve the representation $\rep{p}{0}$ and permute the subalgebras
$W(\NN)$. Namely, $\spfn{p\theta}$ preserves $N_0, \dots, N_{d-1}$,
increases $N_+$ by $p\theta$ and decreases $N_-$ by $p\theta$. 

Recall that $L_p = N_+ + N_- - N_{p-2}$.

First consider the case $L_p >0$. Then we should set
$d=p-1$. Let $\theta \in \oZ$ be the
minimal number such
that $2N_+ + p \theta - N_{p-2} \geq -p+1$.  
So we have $2N_+ + p \theta - N_{p-2}\le p$ and 
as $L_{p}\geq 1$ we obtain $2N_- - p \theta-N_{p-2} > -p+1$. 
Hence after applying 
the automorphism $\spfn{p\theta}$ the  well-balance conditions will be
satisfied.

In the case $L_p = 0$ we set $d=p-2$. Let as above $\theta \in \oZ$ be the
minimal number such
that $2N_+ + p \theta - N_{p-3} \geq -p$.  
So we have $2N_+ + p \theta - N_{p-3} < p$ and 
as $L_{p-1}\geq 0$ we obtain $2N_- - p \theta-N_{p-3} > -p$. 
Hence after applying 
the automorphism $\spfn{p\theta}$ the  well-balance conditions will be
satisfied.  
\end{proof}

\begin{Cor}\label{Cor:iso}
Let $0 \le d \le p-1$.
Let ${\theA}$ be the matrix given by ~\eqref{eq:the-matrix}. Then we have
  \begin{equation}
     \Brep{p}{r\vect u} \cong \brep{p}{r}^d.
  \end{equation}
If $N_d = N_{d+1} = \dots = N_{p-2} = N_+ + N_-$ and 
$\vect N$ satisfies the well-balance condition
$$ 2N_+ - N_{d-1} - 2r \ge -(2p-d-2), \qquad 2N_- - N_{d-1} + 2r \ge
-(2p-d-2)$$
then we have
\begin{equation}
\Brep{p}{r\vect u}[W^{\theA}(\vect N)] \cong \brep{p}{r}^d[W(\vect N)]
\end{equation}
\end{Cor}

So for the character of $\rep{p}{r}$ we have a set of formulas indexed
by $0\le d \le p-1$
$${\rm ch}\, \rep{p}{r} (q,z)  = z^{-\frac{r}{p}}\,q^{\frac{r}{2p}(r-p+2)}\, 
\sum_{\vect{n}\in\oZ^\ell}
  z^{\vect{u}\cdot\vect{n}}\,
  q^{\half\vect{n}{\theA}\cdot\vect{n} + \vect{v}\cdot\vect{n}}
   \prod\limits_{a\in I} \frac1{(q)_{\vect{e}_a\cdot\vect{n}}},$$
where the set of indices is $I=\{+,-,0,1,2,\dots,d-1\}$,
$\vect{u}=(1,-1;0,\dots,0)$, $\vect{v} = (p/2 - r -1) \vect{u}$. 

Now discuss the character formulas for coinvariants.

\begin{Thm}\label{th:char}
Consider the character defined as ${\rm ch}\, 
\rep{p}{r}[W(\NN)](q,z)={\rm tr}\,z^{\cH_0}\,q^{\cL_0}$.

\def \theenumi {\roman{enumi}}
\def\labelenumi {(\theenumi)}
\begin{enumerate}
\item Let the set of indices be $I=\{+,-,0,1,\dots,p-2\}$. Consider the
vectors
$\vect{u}=(1,-1;0,\dots,0)$, $\vect{v} = (p/2 - r -1) \vect{u}$, $\vect{w}
= r\vect{u}$. Then
$${\rm ch}\, \rep{p}{r}[W(\NN)] (q,z) =
z^{-\frac{r}{p}}\,q^{\frac{r}{2p}(r-p+2)}\, 
\sum_{\vect{n}\in\oZ^\ell}
  z^{\vect{u}\cdot\vect{n}}\,
  q^{\half\vect{n}{\theA}\cdot\vect{n} + \vect{v}\cdot\vect{n}}
 \prod\limits_{a\in I}
 \qbin{\vect{e}_a\cdot(\vect{N}+\vect{n}-\vect{n}{\theA}+\vect{w})}
                   {\vect{e}_a\cdot\vect{n}}{q}^+,$$
where the notion $\qbin{n}{m}{q}^+$ is given by~\eqref{eq:qbinplus}.
\item We have
$${\rm ch}\, \rep{p}{r}[W(\NN)] (q,z)=
z^{-\frac{r}{p}}\,q^{\frac{r}{2p}(r-p+2)}\, \sum_{m\in \oZ} z^m
q^{\frac{p}{2}(m^2+m) -(r+1)m}
\qsup{\LL}{pm-r+N_-}{q}^{p+1}.$$
\end{enumerate}
\end{Thm}
\begin{proof}
Note that $\spfn{\theta}$ acts on the
character by substituting $zq$ into $z$. Therefore it is enough to prove
the theorem in the case $r=0$.

First suppose that $\NN$ satisfies the well-balance
conditions~\eqref{wellb}. Then (i) follows from Corollary~\ref{Cor:iso}
and Proposition~\ref{Prop:char}; the statement (ii) can be deduced from
(i) by
Theorem~\ref{ta}. 

We know that general case can be reduced to the well-balanced case by the
action of $\spfn{p\theta}$. As $\spfn{p\theta}$ acts on the RHS of (ii) by
the same substitution $zq^p$ into $z$, we have (ii) in general case. The
statement (i) can be deduced from     
(ii) by
Theorem~\ref{tb}.
\end{proof}

\begin{Cor}
Let $0 \le d \le p-2$. 
Suppose that $N_d = N_{d+1} = \dots = N_{p-2} =
N_+ + N_-$. 
Let the set of indices be $I=\{+,-,0,1,\dots,d-1\}$ and, as above,
$\vect{u}=(1,-1;0,\dots,0)$, $\vect{v} = (p/2 - r -1) \vect{u}$, $\vect{w}
= r\vect{u}$. Then
$${\rm ch}\, \rep{p}{r}[W(\NN)] (q,z) =
z^{-\frac{r}{p}}\,q^{\frac{r}{2p}(r-p+2)}\, 
\sum_{\vect{n}\in\oZ^\ell}
  z^{\vect{u}\cdot\vect{n}}\,
  q^{\half\vect{n}{\theA}\cdot\vect{n} + \vect{v}\cdot\vect{n}}
 \prod\limits_{a\in I}
 \qbin{\vect{e}_a\cdot(\vect{N}+\vect{n}-\vect{n}{\theA}+\vect{w})}
                   {\vect{e}_a\cdot\vect{n}}{q}^+,$$
where the notion $\qbin{n}{m}{q}^+$ is given by~\eqref{eq:qbinplus}.
\end{Cor}

\bigskip

\noindent
\textbf{Acknowledgments}.
This work is supported by grants RFBR 00-15-96579, 01-01-00906, 01-02-16686
and 01-01-00546 and CRDF RP1-2254.

\end{document}